\documentclass[11pt]{article}
	\usepackage{amsmath}
      \usepackage{amssymb}
      \usepackage{pxfonts}
      \usepackage{physics}
        \usepackage{enumitem}
        \usepackage{ragged2e}
        \usepackage{graphicx}
        \usepackage{caption}
        \usepackage{subcaption}
        \usepackage{float}
        \usepackage{mathrsfs}
        \usepackage{mathtools}
        \usepackage{hyperref}
        \usepackage[dvipsnames]{xcolor}
        \usepackage[utf8]{inputenc}
        \usepackage{cite} 
        \usepackage[margin=1in, top=1in, bottom=1in]{geometry}
        \usepackage[a-1b]{pdfx}
        
        \numberwithin{equation}{section}
        \newtheorem{theorem}{Theorem}
        
        \usepackage{lineno}
        \usepackage{authblk}   
\usepackage{hyperref} 
	\begin{document}
	
        \justifying
        \title{Global Analysis of the Gray--Scott Model with Fractional-- Classical Diffusion.}
        \author{Md Shah Alam\thanks{Email: \href{mailto:malam@htu.edu}{malam@htu.edu}}}
\affil{Department of Mathematics, Huston--Tillotson University, Austin, TX 78702, USA}

\date{}
\maketitle{}
\vspace{-0.5 in}
\begin{abstract}
\noindent
We analyze the Gray–Scott reaction–diffusion system on $\Omega \subset \mathbb{R}^n$ $(n \ge 1)$ with mixed diffusion combining local and nonlocal operators. Using semigroup theory and duality estimates, we prove global existence of component-wise nonnegative solutions and establish uniform bounds. Numerical simulations illustrate pattern formation and highlight qualitative differences between the purely local and mixed-diffusion models.
\end{abstract}

\section{Introduction and Main Results}\label{Introduction and Main Results}
\subsection{Introduction}
Classical diffusion processes, governed by integer-order Laplacian operators, are foundational in the modeling of transport phenomena across physics, biology, and engineering \cite{fife2013mathematical,maini1997spatial,enderling2014mathematical,gardiner2004handbook,turing1990chemical,kitanidis1997introduction}. These models rest on the assumption that diffusion is local and governed by Brownian motion, where particles undergo random walks with infinitesimal steps. However, an increasing body of empirical and theoretical evidence reveals that many real-world systems exhibit nonlocal and anomalous diffusion, where standard models fail to capture the observed dynamics.\\
\noindent
Fractional diffusion operators, typically expressed through the fractional Laplacian $(-\Delta)^s$ for $s \in (0,1)$, have emerged as powerful tools to describe such phenomena. Unlike the classical Laplacian, which reflects local interactions, fractional operators encode long-range dependencies by incorporating jump processes, Lévy flights, or heavy-tailed step distributions into the dynamics. These operators arise naturally in the context of continuous-time random walks, kinetic theory, image processing, porous media, and biological pattern formation \cite{metzler2000random,tarasov2011fractional,laskin2002fractional,zaslavsky2002chaos,magin2006fractional,yang2011novel}.\\
\noindent
Mathematically, fractional diffusion operators pose both rich structure and analytical challenges. Their nonlocal nature implies that the evolution at a point depends on the global behavior of the function, requiring new frameworks for function spaces, boundary conditions, and numerical approximation. Moreover, fractional operators exhibit a remarkable interplay between regularization and singularity: they smooth functions like classical elliptic operators, yet they do so in a globally coupled fashion.\\
\noindent
In recent years, the study of partial differential equations (PDEs) involving fractional diffusion has grown into a vibrant area of research. From establishing well-posedness and regularity of solutions to exploring pattern formation, front propagation, and long-time dynamics, fractional PDEs provide a fertile ground for both analysis and applications. In particular, the introduction of fractional Laplacian into reaction-diffusion systems, fluid dynamics, and financial mathematics has led to significant breakthroughs, revealing behaviors absent in classical models.\\
\\
\noindent
The Gray-Scott model is a prototypical reaction--diffusion system that captures autocatalytic chemical reactions with nonlinear feedback, and is defined by a pair of coupled partial differential equations for two species: the activator \( u(x,t) \) and the inhibitor \( v(x,t) \). In its classical form, the model reads:
\begin{flalign} \label{eq:local}
\begin{cases}
    u_t &= d_1 \Delta u - uv^2 + f(1 - u), \\
    v_t & = d_2 \Delta v + uv^2 - (f + k)v,
    \end{cases}
\end{flalign}
where \( \Delta \) denotes the standard Laplacian operator acting on a bounded spatial domain \( \Omega \subset \mathbb{R}^n \), typically subject to homogeneous Neumann or periodic boundary conditions. The parameters \( d_1, d_2 > 0 \) are diffusion coefficients, while \( f \) and \( k \) represent feed and removal rates, respectively. The model originates from the cubic autocatalysis scheme
\begin{flalign*}
     U + 2V & \to 3V\\
     V & \to P
\end{flalign*}
introduced by Gray and Scott (1984), and has since become a canonical example in the study of pattern formation, due to its ability to generate a wide range of spatial structures including spots, stripes, and labyrinthine patterns. The Laplacian term models Fickian diffusion, ensuring spatial coupling between neighboring concentrations. Mathematically, the system exhibits rich dynamics and has been extensively studied for its Turing instability, global existence of solutions, and numerical pattern formation, particularly in two-dimensional domains. Its accessibility and visual expressiveness have also made it a frequent testbed in numerical analysis and computational PDE research.\\
\noindent
In parallel, there has been substantial progress on \emph{nonlocal} and \emph{fractional} diffusion operators, both in terms of the correct functional setting and the appropriate notion of boundary condition. In particular, fractional Laplacians with nonlocal Neumann-type conditions (formulated via the corresponding energy form and nonlocal flux) have been studied systematically; see \cite{dipierro2017nonlocal,cinti2020nonlocal}. 
Regarding Gray--Scott-type dynamics with nonlocal diffusion, well-posedness and nonlocal-to-local limits (diffusive limits) have recently been established for convolution-based nonlocal Gray--Scott models; see \cite{laurencot2023nonlocal}. Moreover, fractional diffusion is known to significantly affect pattern selection and Turing-type mechanisms in reaction--diffusion settings; \cite{langlands2007turing,gafiychuk2006pattern}.\\
\\
\noindent
Motivated by applications where different species may diffuse through different physical mechanisms, we study a \emph{mixed} Gray--Scott system in which the activator $u$ diffuses by a fractional Laplacian $( -\Delta)^s$ while the inhibitor $v$ diffuses by the classical Laplacian, both posed on a smooth bounded domain and complemented with homogeneous Neumann-type conditions (nonlocal Neumann for $u$ and classical Neumann for $v$). 
The main novelty is that the model couples \emph{two different diffusion generators} (fractional/nonlocal versus classical/local) in a single reaction--diffusion system and we provide a complete global well-posedness and uniform boundedness theory in this mixed setting.\\
\noindent
Our approach is operator-theoretic and proceeds in two steps. First, we realize the fractional diffusion operator through its closed, symmetric bilinear form and the associated self-adjoint operator, while the Neumann Laplacian is treated by standard elliptic theory. This yields analytic semigroups for each component, allowing the system to be written as a coupled semilinear evolution problem on the product space $L^2(\Omega)\times L^2(\Omega)$ with locally Lipschitz reaction terms. Second, to rule out finite-time blow-up in the $L^\infty$-norm, we combine (i) a semigroup comparison estimate to obtain an explicit uniform bound for $u$ and an $L^1$-bound for the total density $u+v$, with (ii) a duality bootstrap argument to upgrade integrability and ultimately derive uniform-in-time $L^\infty$ bounds for $v$ as well.\\
\noindent
Finally, we complement the analysis with numerical experiments that compare the classical Gray--Scott patterns with those produced by the mixed diffusion model and illustrate how the fractional order $s$ changes morphology, with patterns approaching the classical ones as $s\to 1^{-}$ (as observed in Figures~1--4).\\
\\
\noindent
Let $n \in \mathbb{N}$ and $\Omega$ is a bounded open subset of $\mathbb{R}^n$ with smooth boundary $\partial \Omega$ ($\partial \Omega$ is an $n-1$ dimensional $C^{2+\mu}$ $(\mu \in (0,1))$ manifold such that $\Omega$ lies locally on one side of $\partial \Omega$). We consider the following Gray-Scott model with both local and nonlocal diffusion operator
\begin{flalign}\label{eq:1.1}
\begin{cases}
u_t(x,t) & = - d_1 (-\Delta)^{s} u(x,t) - u(x,t) v(x,t)^2 + f(1- u(x,t)), \ \ \ (x,t) \in \Omega \times (0,\infty),\\
v_t(x,t) & = d_2 \Delta v(x,t) + u(x,t) v(x,t)^2 - (f+ \kappa)v(x,t), \hspace{0.37 in} (x,t) \in \Omega \times (0,\infty),\\
\mathcal{N}_s u & = 0, \hspace{2.84 in} (x,t) \in \mathbb{R}^n \setminus \Omega \times (0,\infty),\\
\frac{\partial v}{\partial \eta} & = 0, \hspace{2.84 in} (x,t) \in \partial \Omega \times (0,\infty),\\
u(x,0) & = u_0(x) \ge 0, \ v(x,0) =v_0(x) \ge 0\hspace{1.16  in} x \in \Omega,
\end{cases}
\end{flalign}
\noindent
$u$ and $v$ are concentrations of two chemicals with nonnegative initial data $u_0,v_0 \in L^2(\Omega)$. $(-\Delta)^s$ is nonlocal fractional Laplacian with $s \in (0,1)$ associated with nonlocal homogeneous Neumann boundary condition,
\begin{flalign*}
    \mathcal{N}_s u = C_{n,s} \int_{\Omega} \frac{u(y,t)-u(x,t)}{\vert x-y \vert^{n+2s}} \ dy, \ \ \ (x,y) \in \mathbb{R}^n \setminus \Omega \times \mathbb{R}^n \setminus \Omega, \ x \neq y, \ t>0
\end{flalign*}
in the Cauchy Principle Value sense so that the integral does not diverge. $d_1,d_2>0$ are diffusion coefficients with $d_1 \neq d_2$, $f$ is called feed rate and $\kappa$ is called the kill rate of $u$. Above model consists of two different types of diffusion operators, a faster diffusion operator (\textit{fractional Laplacian}) and a slower (\textit{regular/local Laplacian}) one. The diffusion of the fractional one is anomalous in nature while the other one is regular. Let us consider the following locally Lipschitz functions
\begin{flalign*}
g_1(u,v) & = - uv^2 + f(1- u)\ \ \ \text{and } \ g_2(u,v) = v^2 - (f+ \kappa)v
\end{flalign*}
with $g_1(0,v),g_2(u,0) \ge 0$ whenever $u,v \ge 0$. This condition is known as \textit{quasi-positivity} and it ensures the nonnegativity of the solution along with nonnegative initial data.\\
\\
\noindent
Let $\mathbb{R}^2_+=[0,\infty)^2$. We consider the following Riesz definition of fractional Laplacian in the Cauchy Principle Value sense so that the integral does not diverge
\begin{flalign}\label{eq:1.2}
    (-\Delta)^s u (x,t) & = C_{n,s} \int_{\Omega} \frac{u(x,t)-u(y,t)}{\vert x-y \vert^{n+2s}} \ dy, \ \ \ (x,y) \in \Omega \times \Omega, \ x \neq y, \ t>0
\end{flalign}
where $C_{n,s} = \frac{4^s \Gamma (\frac{n}{2}+s)}{\pi^{\frac{n}{2}}\Gamma(-s)}$ with $\Gamma(z) = \int_0^{\infty} t^{z-1} e^{-t} \ dt , \ Re(z)>0$.\\
\noindent
Throughout this work, the fractional diffusion operator acting on the activator $u$ is understood as the \emph{Neumann realization of the fractional Laplacian} on $\Omega$.\\
\noindent
More precisely, for $s\in(0,1)$, the fractional Laplacian $(-\Delta)^s$ is defined via the associated bilinear form
\begin{equation*}\label{eq:bilinear_form}
\mathcal{E}_s(u,\varphi)
=
\frac{C_{n,s}}{2}
\int_{\Omega}\int_{\Omega}
\frac{(u(x)-u(y))(\varphi(x)-\varphi(y))}{|x-y|^{n+2s}}
\,dx\,dy,
\end{equation*}
for all admissible test functions $\varphi\in H^s(\Omega)$.
The normalization constant $C_{n,s}>0$ is chosen such that, when $\Omega=\mathbb{R}^n$, the operator coincides with the classical fractional Laplacian defined through the Fourier multiplier $|\xi|^{2s}$.\\
\noindent
The operator $(-\Delta)^s$ is then characterized as the unique self-adjoint operator on $L^2(\Omega)$ associated with the closed, coercive bilinear form $\mathcal{E}_s$, and its domain is given by
\begin{equation*}
D((-\Delta)^s)
=
\left\{
u\in H^s(\Omega):
\; (-\Delta)^s u \in L^2(\Omega)
\right\}.
\end{equation*}
\noindent
The \emph{nonlocal Neumann boundary condition} is incorporated intrinsically through the variational formulation of $\mathcal{E}_s$. Formally, it corresponds to the vanishing of the nonlocal normal flux across the boundary,
\begin{equation*}\label{eq:nonlocal_neumann}
\mathcal{N}_s u(x)
=
C_{n,s}
\int_{\Omega}
\frac{u(y)-u(x)}{|x-y|^{n+2s}}
\,dy
=0,
\qquad x\in\partial\Omega,
\end{equation*}
which ensures conservation of total mass in the absence of reaction terms.
\noindent
With this realization, the operator $(-\Delta)^s$ is nonnegative, self-adjoint, and generates an analytic semigroup on $L^2(\Omega)$. This framework is consistent with the analytical results developed in Section~2 and with the numerical implementation described in Section~3. Further details on nonlocal Neumann problems and related fractional operators can be found in \cite{cinti2020nonlocal,dipierro2017nonlocal}.

\subsection{Statement of the Main Results}
\begin{theorem} \label{Th:1} If $(u_0,v_0) \in L^2(\Omega,\mathbb{R}^2_+)$ with $\frac{\partial}{\partial \eta} v_0=0$ on $\partial \Omega=0$ then there exists a unique classical global componentwise nonnegative solution $(u(x,t),v(x,t))$ to (\ref{eq:1.1}) and both $u$ and $v$ remain uniformly bounded in the sup norm for all time.
\end{theorem}
\noindent
In section (\ref{Proof of Main Results}), we prove our main  theorem. Then, we present our numerical illustration of patterns in section (\ref{Pattern Formation}). Conclusions are discussed in section (\ref{Conclusion and Future Directions}).

\section{Proof of Main Results}\label{Proof of Main Results}
Define $X := L^2(\Omega)$, $X_+ := \{ z \in X: z \ge 0 \text{ in } \Omega\}$. For the sake of simplicity of notation, we also define $A_1 (z)=d_1(-\Delta)^s z$, $A_2 (z)= -d_2 \Delta z$ with the following domain respectively for all
\begin{flalign*}
 z \in D(A_1) & = \big\{w\in H^s(\Omega)\,: w \in L^2(\Omega), \ A_1(w) \in L^2(\Omega),\,\mathcal{N}_s u=0\,\text{on }\mathbb{R}^n \setminus \Omega \big\}\\
z\in D(A_2) & = \big\{w\in W^{2,2}(\Omega)\,: \ A_2(w)\in L^2(\Omega),\,\frac{\partial w}{\partial\eta}=0\,\text{on }\partial\Omega \big\},
\end{flalign*}
where the fractional Sobolev space is defined as
\begin{flalign*}
   H^s(\Omega) := \big\{w \in L^2(\Omega): \int_{\Omega} \int_{\Omega} \frac{\vert u(y)-u(x) \vert^2}{\vert x-y \vert^{n+2s}} dxdy<\infty  \big\} 
\end{flalign*}
with norm
\begin{flalign}\label{eq:H_norm}
    \Vert u \Vert^2_{H^s(\Omega)} := \Vert u \Vert^2_{L^2(\Omega)}+\int_{\Omega} \int_{\Omega} \frac{\vert u(y)-u(x) \vert^2}{\vert x-y \vert^{n+2s}} = \Vert u \Vert^2_{L^2(\Omega)} + [u]^2_{H^s(\Omega)}.
\end{flalign}
$[u]_{H^s(\Omega)}$ is a semi-norm introduced by Gagliardo.\\
\noindent
The operator $A_2$ is well known to be self-adjoint, nonnegative, and densely defined on $L^2(\Omega)$ with domain $D(A_2)$. Moreover, $A_2$ is sectorial and therefore generates an analytic semigroup $\{e^{tA_2}\}_{t\ge 0}$ on $L^2(\Omega)$\cite{hollis1987global}.\\
\noindent
For the fractional operator $A_1$, we consider the bilinear form
\[
\mathcal{E}_s(u,\varphi)
=
\frac{d_1 C_{n,s}}{2}
\int_{\Omega}\int_{\Omega}
\frac{(u(x)-u(y))(\varphi(x)-\varphi(y))}{|x-y|^{n+2s}}
\,dx\,dy,
\]
defined on $H^s(\Omega)\times H^s(\Omega)$.
Then this form is symmetric, closed, and coercive on $L^2(\Omega)$. Since $\mathcal{E}_s$ is a densely defined, closed, symmetric, and nonnegative bilinear form on $L^2(\Omega)$ with domain $H^s(\Omega)$, the Kato--Friedrichs representation theorem yields a unique self-adjoint, nonnegative operator $A_1$ associated with $\mathcal{E}_s$. In particular, $A_1$ is sectorial and generates an analytic semigroup $\{e^{tA_1}\}_{t\ge 0}$ on $L^2(\Omega)$ \cite{daoud2024class}.\\
\noindent
Now, \ref{eq:1.2} is therefore naturally formulated as a coupled semilinear evolution equation on the product space $L^2(\Omega)\times L^2(\Omega)$. Since, each component evolves under its own diffusion semigroup and the coupling occurs exclusively through the locally Lipschitz reaction terms, therefore different boundary conditions, does not have any impact on the semigroup generated by the coupled system.
Defining the block-diagonal operator
\[
\mathcal{A}
=
\begin{pmatrix}
A_1 & 0 \\
0 & A_2
\end{pmatrix},
\qquad
D(\mathcal{A})=D(A_1)\times D(A_2),
\]
standard semigroup theory implies that $\mathcal{A}$ generates a strongly continuous (and analytic) product semigroup
\[
e^{t\mathcal{A}}(u_0,v_0)
=
\bigl(e^{tA_1}u_0,\; e^{tA_2}v_0\bigr)
\quad \text{on } L^2(\Omega) \times L^2(\Omega),
\]
see, for instance, \cite{pazy2012semigroups,arendt1987vector}
.\\
\\
\noindent
\textbf{Proof of Theorem {\ref{Th:1}}:} With the above settings, results in \cite{pazy2012semigroups} guarantee that (\ref{eq:1.1}) has a unique maximal solution $(u,v)$ on a maximal interval $[0,T_m)$, where $T_m$ is a positive extended real number, given by 
\begin{flalign*}
    u(x,t) & = T_1(t) u_0(x) +\int_0^t T_1(t-\tau) g_1(u(x,\tau)v(x,\tau)) d\tau \\
    v(x,t) & = T_2(t) v_0(x) +\int_0^t T_2(t-\tau) g_2(u(x,\tau),v(x,\tau)) d\tau
\end{flalign*}
It is straightforward to prove that \textit{quasi-positivity} implies the maximal solution $(u,v)$ is component-wise nonnegative. To establish the global existence, we need to show either $T_m=\infty$ or 
\begin{flalign}
  \limsup_{t \to T_m^-} \big\{ \Vert u(\cdot,t)\Vert_{\infty,\Omega} + \Vert v(\cdot,t)\Vert_{\infty,\Omega} \big\} = \infty. \label{eq:2.1}
\end{flalign}
Consequently, we demonstrate global existence by proving solutions cannot blow up in the sup norm in finite time.
Now, if we define $A_{11}=A_1+fI=d_1(-\Delta)^s+fI$, then $A_{11}(u)$ generates $T_{11}(t)$, and from results in \cite{laurencot2023nonlocal} we observe that $\Vert T_{11}(z) \Vert_{\infty,\Omega} \le e^{-tf}\Vert z \Vert_{\infty,\Omega}$ for all $z \in L^{\infty}(\Omega,\mathbb{R})$. In addition,
\begin{flalign}
u_t(x,t) & = - d_1(-\Delta)^s u(x,t) - u(x,t) v(x,t)^2 + f(1- u(x,t)) \nonumber\\
=> \ u(x,t) & \le T_{11}(t) u_0(x) + \int_0^t T_{11}(t-\tau) f \ d\tau \hspace{1.4 in} [\because \ u,v\ge 0 \ \forall \ t>0] \nonumber\\
=> \ u(x,t) & \le \Vert T_{11}(t) u_0 \Vert_{\infty,\Omega} + \int_0^t \Vert T_{11}(t-\tau) f \Vert_{\infty,\Omega} \ d\tau \le e^{-tf} \Vert u_0 \Vert_{\infty,\Omega} + 1 - e^{-tf} \nonumber\\
\therefore \ \Vert u(\cdot,t) \Vert_{\infty,\Omega} & \le C=\max\big\{\Vert u_0 \Vert_{\infty,\Omega}, 1 \big\}, \ \ \forall t>0  \label{eq:2.2}
\end{flalign}
which is uniform in time and depends only on the initial data $u_0$ and the parameter $f$.\\
\noindent
Now, adding the first two equations of the (\ref{eq:1.1}) we get,
\begin{flalign*}
  u_t(x,t) + v_t(x,t) & = -d_1 (-\Delta)^s u(x,t) + d_2 \Delta v(x,t) + f(1- u(x,t)) -(f+\kappa) v(x,t)\\
=> \ \frac{\partial}{\partial t} \big(u(x,t) + v(x,t)\big) & \le -d_1 (-\Delta)^s u(x,t) + d_2 \Delta v(x,t) + f - f\tilde{\kappa}(u(x,t)+v(x,t))\\
=> \ \frac{d}{dt} \int_{\Omega} \big(u(x,t) + v(x,t)\big) dx & \le -d_1 \int_{\Omega} (-\Delta)^s u(x,t) dx + d_2 \int_{\Omega} \Delta v(x,t) dx + f \int_{\Omega} dx \\
& - f\tilde{\kappa} \int_{\Omega} (u(x,t)+v(x,t)) dx\\
=> \ \frac{d}{dt} \int_{\Omega} \big(u(x,t) + v(x,t)\big) dx & + f\tilde{\kappa} \int_{\Omega} \big(u(x,t)+v(x,t) \big) dx \le  f \vert \Omega \vert 
\end{flalign*}
where $\tilde{\kappa}=\min\big\{ \kappa,1 \big\}>0$ and we have used boundary conditions with Green's formula to get
\begin{flalign*}
d_1 \int_{\Omega} (-\Delta)^s u(x,t) dx = -\int_{\mathbb{R}^n \setminus \Omega} \mathcal{N}_s u(x,t) =0 \ \ \text{and } \ d_2 \int_{\Omega} \Delta v(x,t) dx =0.
\end{flalign*}
Solving the above differential inequality we get,
\begin{flalign}
\int_{\Omega} \big(u(x,t) + v(x,t)\big) dx & \le M\ \ \ \forall t>0. \ \ \ \ \ \big[M = \max\big\{\frac{\vert \Omega \vert}{\tilde{\kappa}},\int_{\Omega}(u_0(x)+v_0(x))dx \big\} \big]\label{eq:2.3}
\end{flalign}
By (\ref{eq:2.2}) and the boundedness of $\Omega$, we have, for every $1\le p<\infty$ and $0\le \tau<T<T_m$,
\begin{flalign}\label{eq:u_p}
\|u\|_{p,\Omega\times(\tau,T)} & =
\left(\int_{\tau}^{T}\int_{\Omega}|u(x,t)|^p\,dxdt\right)^{\frac1p} \le |\Omega|^{\frac1p}(T-\tau)^{\frac1p} \|u\|_{\infty,\Omega\times(\tau,T)} \le M_p(T-\tau),
\end{flalign}
where $M_p(T-\tau)>0$ is independent of $u$ and $t$.

\noindent
We will use (\ref{eq:u_p}) to bootstrap the $L^1$ bounds in (\ref{eq:2.3}) to better estimates ($L^p(\Omega \times (0,T_m))$ for all $1<p<\infty$) by applying duality arguments. To this end, let $0 \le \tau<T<T_m$, with $\tau$ further specified below, $1<q<\infty$ and $\theta \in L^q(\Omega \times (\tau,T))$ such that $\theta \ge0$ and $\|\theta\|_{q,\Omega\times(\tau,T)}=1$. Let $\varphi$  be the unique nonnegative solution in $W_q^{(2,1)}(\Omega\times(\tau,T))$ solving the following system on the left, which at a first glance, appears to be a backward heat equation. But the substitutions $\psi(x,t)=\varphi(x,T+\tau-t)$ and $\tilde\theta(x,t)=\theta(x,T+\tau-t)$ lead to the standard initial boundary value problem given in the right 
\begin{flalign} \label{eq:2.4}
\begin{aligned}[c]
 \begin{cases}
    \varphi_t & = - d_2 \Delta \varphi-\theta, \ \ \ \ \ \ \Omega \times (\tau,T),\\
    \frac{\partial \varphi}{\partial \eta} & = 0, \hspace{0.77 in} \partial\Omega \times (\tau,T),\\
    \varphi & = 0, \hspace{0.85 in} \Omega \times \{T\}.
\end{cases}   
\end{aligned} \ \ \ \ \ \ \ 
\begin{aligned}[c]
\begin{cases}
    \psi_t & = d_2 \Delta \psi+\tilde\theta(x,T+\tau-t), \ \ \ \ \ \ \Omega \times (\tau,T),\\
    \frac{\partial \psi}{\partial \eta} & = 0, \hspace{1.48 in} \partial\Omega \times (\tau,T),\\
    \psi & = 0, \hspace{1.56 in} \Omega \times \{\tau\}.
\end{cases}   
\end{aligned}
\end{flalign}
The duality argument is used to upgrade the $L^1$ bounds on the total density $u+v$ to higher integrability.
Given $\theta\in L^q(\Omega\times(\tau,T))$ with $\|\theta\|_{L^q}=1$, we introduce the auxiliary function $\varphi$ as the solution of the backward parabolic problem
\[
-\partial_t \varphi - d_2 \Delta \varphi = \theta
\quad \text{in } \Omega\times(\tau,T),
\]
supplemented with homogeneous Neumann boundary conditions and terminal condition $\varphi(\cdot,T)=0$.
This choice allows us to transfer derivatives from $u+v$ onto $\varphi$ via integration by parts, so that the forcing term $\theta$ appears explicitly in the resulting estimate.\\
\noindent
Results in \cite{ladyzhenskaia1968linear} imply there exists $C_{q,T-\tau}>0$ so that $\|\varphi\|_{q,\Omega\times(\tau,T)}^{(2,1)}\le C_{q,T-\tau}$ and Lemma 4.1 in \cite{morgan1989global} gives $\|\varphi(\tau)\|_{q,\Omega}\le C_{q,T-\tau}$. Now, we start with
\begin{flalign} 
 & \int_{\tau}^T \int_{\Omega} \big(u(x,t)+v(x,t) \big) \theta(x,t) \ dxdt = \int_{\tau}^T \int_{\Omega} \big(u(x,t)+v(x,t)\big) \big(-d_2\Delta \varphi (x,t)-\varphi_t(x,t)\big) \ dxdt \nonumber\\
  & = -d_2 \int_{\tau}^T \int_{\Omega}  \big(u(x,t) \Delta \varphi(x,t) + v(x,t) \Delta \varphi(x,t) \big) \ dxdt - \int_{\tau}^T \int_{\Omega} \big(u(x,t)\varphi_t(x,t) + v(x,t) \varphi_t(x,t) \big) \ dxdt \nonumber\\
& = A+B \label{eq:2.5}
\end{flalign}
where,
\begin{flalign*}
A & = -d_2 \int_{\tau}^T \int_{\Omega}  (u(x,t) \Delta \varphi(x,t) + v(x,t) \Delta \varphi(x,t)) \ dxdt \\
 & = -d_2 \int_{\tau}^T \int_{\Omega}  u(x,t) \Delta \varphi(x,t) \ dxdt -d_2 \int_{\tau}^T \int_{\Omega} v(x,t) \Delta \varphi(x,t) \ dxdt \\
 & \le d_2 \int_{\tau}^T \int_{\Omega} \vert u(x,t) \Delta \varphi(x,t) \vert \ dxdt -d_2 \int_{\tau}^T \int_{\Omega} v(x,t) \Delta \varphi(x,t) \ dxdt \\
 & \le d_2 \Vert u(\cdot,t) \Vert_{p,\Omega \times (\tau,T)} C_{q,T-\tau} -d_2 \int_{\tau}^T \int_{\Omega} v(x,t) \Delta \varphi(x,t) \ dxdt \\
\therefore \ A & \le d_2 M_p(T-\tau) C_{q,T-\tau} - d_2 \int_{\tau}^T \int_{\Omega} v(x,t) \Delta \varphi(x,t) \ dxdt
\end{flalign*}
and
\begin{flalign*}
B & = - \int_{\tau}^T \int_{\Omega}  (u(x,t) \varphi_t(x,t) + v(x,t) \varphi_t(x,t)) \ dxdt \nonumber\\
 & = - \int_{\Omega} \big[\big(u(x,t)+v(x,t) \big)\varphi(x,t) \big]_{\tau}^T \ dx + \int_{\tau}^T \int_{\Omega} \big(u_t(x,t)+v_t(x,t) \big) \varphi(x,t) \ dxdt \nonumber\\
 & = - \int_{\Omega} \big[\big(u(x,T)+v(x,T)\big)\varphi(x,T)-\big(u(x,\tau)+v(x,\tau)\big)\varphi(x,\tau) \big] \ dx + \int_{\tau}^T \int_{\Omega} \big[\big(-d_1 (-\Delta)^s u(x,t) \\
 & + d_2 \Delta v(x,t)\big)-f (u(x,t)+v(x,t))+ f-\kappa v(x,t)\big] \varphi(x,t) \ dx \ dt \nonumber\\
 & \le \int_{\Omega} \big(u(x,\tau)+v(x,\tau)\big)\varphi(x,\tau) \ dx + \int_{\tau}^T \int_{\Omega} \big(-d_1 (-\Delta)^s u(x,t)  + d_2 \Delta v(x,t)\big) \varphi(x,t) \ dxdt\\
 & + f \int_{\tau}^T \int_{\Omega} \varphi(x,t) \ dxdt - f\tilde{\kappa}\int_{\tau}^T \int_{\Omega}  (u(x,t)+v(x,t)) \varphi(x,t) \ dxdt\\
 & \le \int_{\Omega} \big(u(x,\tau)+v(x,\tau)\big)\varphi(x,\tau) \ dx - \int_{\tau}^T \int_{\Omega} d_1 (-\Delta)^s u(x,t) \varphi(x,t) \ dxdt\\
 & + \int_{\tau}^T \int_{\Omega} d_2 \Delta v(x,t)\varphi(x,t) \ dxdt + f \int_{\tau}^T \int_{\Omega} \varphi(x,t) \ dxdt - f\tilde{\kappa}\int_{\tau}^T \int_{\Omega}  (u(x,t)+v(x,t)) \varphi(x,t) \ dxdt
\end{flalign*}
If $q=\frac{n+4}{2}$ then $W^{(2,1)}_q(\Omega \times (\tau, T))$ embeds continuously into $C(\overline\Omega \times [\tau, T])$ \cite{ladyzhenskaia1968linear}.\\
\noindent
By standard parabolic regularity theory, $\varphi\in W^{(2,1)}_q(\Omega\times(\tau,T))$, with a bound depending only on $\|\theta\|_{L^q}$.
The exponent $q=\frac{n+4}{2}$ is chosen so that the parabolic Sobolev embedding $W^{(2,1)}_q(\Omega\times(\tau,T)) \hookrightarrow C(\overline{\Omega}\times[\tau,T])$ holds \cite{ladyzhenskaia1968linear}.
This embedding guarantees uniform boundedness of $\varphi$ and of its time trace, which is essential for controlling reaction terms and initial-time contributions in the duality estimate.
Consequently, the pairing $\int_{\tau}^{T}\!\!\int_{\Omega} (u+v)\,\theta$ is uniformly bounded, yielding higher integrability of $u+v$ by duality.\\
\noindent
Using this estimate along with (\ref{eq:2.3}) we get,
\begin{flalign*}
 \int_{\Omega} \big(u(x,\tau)+v(x,\tau)\big) \varphi(x,\tau) dx & \le K_{q,T-\tau} \int_{\Omega} \big(u(x,\tau) + v(x,\tau)\big) dx \le K_{q,T-\tau}M, \\
 \int_{\tau}^T \int_{\Omega} \big(u(x,t)+v(x,t)\big) \varphi(x,t) dxdt & \le K_{q,T-\tau} M(T-\tau) \ \ \text{and } \ \int_{\tau}^T \int_{\Omega} \varphi(x,t) dxdt \le M_{q,T-\tau}
\end{flalign*}
where for the last integral, we have used the Trace Embedding Theorem in \cite{ladyzhenskaia1968linear} to get that there exists $M_{q,T-\tau}>0$ independent of $\theta$ such that $\Vert \varphi \Vert_{L^1,\Omega \times (\tau,T)} \le M_{q,T-\tau}$. We also observe that,
\begin{flalign*}
\int_{\tau}^T \int_{\Omega} d_2 \Delta v(x,t) \varphi(x,t) \ dxdt  = d_2 \int_{\tau}^T \int_{\Omega} v(x,t) \Delta \varphi(x,t) \ dxdt.
\end{flalign*}
Finally,
\begin{flalign}
 & - \int_{\tau}^T \int_{\Omega} d_1 (-\Delta)^{s} u(x,t) \varphi(x,t) \ dxdt\nonumber\\
 & =  - \frac{d_1 C_{n,s}}{2} \int_{\tau}^T \int_{\Omega} \int_{\Omega}  \frac{(u(y,t) - u(x,t))(\varphi(y,t)-\varphi(x,t))}{\vert x-y \vert^{n+2s}} dxdydt + \int_{\mathbb{R}^n \setminus \Omega} \varphi(x,t) \mathcal{N}_{s} u(x,t) dx \nonumber\\
 & \le \frac{d_1 C_{n,s}}{2} \big[\int_{\tau}^T \int_{\Omega} \int_{\Omega}  \frac{\vert u(y,t) - u(x,t) \vert^2}{\vert x-y \vert^{n+2s}} dxdydt \big]^{\frac{1}{2}} \big[\int_{\tau}^T \int_{\Omega} \int_{\Omega}  \frac{\vert\varphi(y,t)-\varphi(x,t) \vert^2}{\vert x-y \vert^{n+2s}}\big]^{\frac{1}{2}} dxdydt + 0 \nonumber\\  
 & \hspace{3.9 in}[\text{applying H\"older's inequality}] \nonumber\\
 & \le \frac{d_1 C_{n,s}}{4} \big[\int_{\tau}^T \int_{\Omega} \int_{\Omega}  \frac{\vert u(y,t) - u(x,t) \vert^2}{\vert x-y \vert^{n+2s}} dxdydt + \int_{\tau}^T \int_{\Omega} \int_{\Omega}  \frac{\vert\varphi(y,t)-\varphi(x,t) \vert^2}{\vert x-y \vert^{n+2s}} dxdydt \big] \label{eq:frac_u}
\end{flalign} 
where we used Young's inequality.

Now, let us multiply both sides of the first equation of system (\ref{eq:1.1}) with $u(x,t)$ and then integrating over $\Omega$ to get,
\begin{flalign*}
\frac{1}{2} \frac{d}{dt} \int_{\Omega} u(x,t)^2 dx & \le -d_1 \int_{\Omega} u(x,t) (-\Delta)^s u(x,t) dx + f \int_{\Omega}  u(x,t) dx - f \int_{\Omega} u(x,t)^2 dx
\end{flalign*}
Also, if we define $\phi(x,y)=\frac{1}{\vert x-y \vert^{n+2s}}$ then using (\ref{eq:1.2}) we can write,
\begin{flalign*}
  - d_1 \int_{\Omega} u(x,t) (-\Delta)^{s} u(x,t) dx = d_1 C_{n,s} \int_{\Omega} \int_{\Omega} u(x,t) \phi(x,y) (u(y,t) - u(x,t)) dydx
\end{flalign*}
since, $\phi(x,y)=\phi(y,x)$ and $\phi(x,y) > 0$ for $(x,y) \in (\Omega \times \Omega)$ with $x \neq y$, we get,
\begin{flalign*}
  - d_1 \int_{\Omega} u(x,t) (-\Delta)^{s} u(x,t) dx =  -\frac{d_1 C_{n,s}}{2} \int_{\Omega} \int_{\Omega} \phi(x,y)  (u(y,t) - u(x,t))^2 dydx.
\end{flalign*}
Also, using Young's inequality, for every $0<\varepsilon<1$, we can find a constant $C_{\varepsilon}>0$ so that
\begin{align*}
f\int_{\Omega}u(x,t)\,dx &\leq \varepsilon f\int_{\Omega}u(x,t)^2\,dx + C_{\varepsilon}f|\Omega|.
\end{align*}
Replacing these estimates we get,
\begin{flalign*}
\frac{1}{2} \frac{d}{dt} \int_{\Omega} u(x,t)^2 dx + \frac{d_1 C_{n,s}}{2} \int_{\Omega} \int_{\Omega} \phi(x,y) (u(y,t) - u(x,t))^2 dydx + f(1-\varepsilon) \int_{\Omega} u(x,t)^2 dx \le f C_{\varepsilon,2}\vert \Omega \vert
\end{flalign*}
for any $0<\varepsilon<1$, integrating the above inequality from $\tau$ to $T$ we can write,
\begin{flalign*}
\frac{d_1 C_{n,s}}{2} \int_{\tau}^T \int_{\Omega} \int_{\Omega} \phi(x,y) (u(y,t) - u(x,t))^2 dydxdt & \le f C_{\varepsilon,2}\vert \Omega \vert (T-\tau)\\
\therefore \ \frac{d_1 C_{n,s}}{4} \int_{\tau}^T \int_{\Omega} \int_{\Omega} \frac{\vert u(y,t) - u(x,t) \vert^2}{\vert x-y \vert^{n+2s}} dydx & \le \frac{f C_{\varepsilon,2}\vert \Omega \vert (T-\tau)}{2}.
\end{flalign*}
Also, using (\ref{eq:H_norm}) along with the the continuous embedding of $W^{2,1}_q \hookrightarrow W^{2,1}_s=H^s$ for $2-\frac{n}{q}>s-\frac{n}{2}$ for any $C_2>0$, we can write,
\begin{flalign*}
\int_{\tau}^T \int_{\Omega} \int_{\Omega}  \frac{\vert\varphi(y,t)-\varphi(x,t) \vert^2}{\vert x-y \vert^{n+2s}} dxdydt & = [\varphi(\cdot,t)]^2_{H^{\alpha}(\Omega \times (\tau,T))} = \Vert \varphi(\cdot,t) \Vert^2_{H^{\alpha}(\Omega \times (\tau,T))} - \Vert \varphi(\cdot,t) \Vert^2_{L^2(\Omega \times (\tau,T))}\\
& \le \Vert \varphi(\cdot,t) \Vert^2_{H^{\alpha}(\Omega \times (\tau,T))} \le C_2 \Vert \varphi(\cdot,t) \Vert^2_{W^{2,1}_q(\Omega \times (\tau,T))} \le C_2 C^2_{q,T-\tau}.
\end{flalign*}
Replacing these estimates on (\ref{eq:frac_u}), we get
\begin{flalign*}
-\int_{\tau}^T \int_{\Omega} d_1 (-\Delta)^{s} u(x,t) \varphi(x,t) \ dxdt & \le \frac{f C_{\varepsilon,2}\vert \Omega \vert (T-\tau)}{2} + \frac{d_1 C_{n,s} C_2 C^2_{q,T-\tau}}{4} 
\end{flalign*}
Replacing these estimates in $B$ and then replacing $A,B$ in  $(\ref{eq:2.5})$ we get,
\begin{flalign*} 
\int_{\tau}^T \int_{\Omega} \big(u(x,t) + v(x,t)\big) \theta(x,t) \ dxdt & \le d_2 M_p(T-\tau) C_{q,T-\tau} + K_{q,T-\tau}M + \frac{f C_{\varepsilon,2}\vert \Omega \vert (T-\tau)}{2} + \frac{d_1 C_{n,s} C_2 C^2_{q,T-\tau}}{4}\\
& + f M_{q,T-\tau} +  f\tilde{\kappa}K_{q,T-\tau}M(T-\tau) \\
& \le N_{q,T-\tau} (1+M)
\end{flalign*}
where $N_{q,T-\tau}= \max \big\{d_2 M_p(T-\tau) C_{q,T-\tau} + \frac{f C_{\varepsilon,2}\vert \Omega \vert (T-\tau)}{2} + \frac{d_1 C_{n,s} C_2 C^2_{q,T-\tau}}{4} + f M_{q,T-\tau} +  f\tilde{\kappa}K_{q,T-\tau}M(T-\tau), K_{q,T-\tau} \big\}$. Now, by choosing $\tau$ such that $T-\tau$ is sufficiently small and using $q'=\frac{n+4}{n+2}$, we find there exists,
\begin{flalign*} 
\Vert v(\cdot,t) \Vert_{\frac{n+4}{n+2},\Omega \times (\tau,T)} \le N_{q,T-\tau} (1+M)
\end{flalign*}
Note that since this estimate only depends on the size of $T-\tau$, it follows from a limiting process that
\begin{flalign*} 
\Vert v(\cdot,t) \Vert_{\frac{n+4}{n+2},\Omega \times (\tau,T_m)} \le N_{q,T_m-\tau} (1+M)
\end{flalign*}
Now, let us write the above estimate as
\begin{flalign} \label{eq:2.6}
\Vert v(\cdot,t) \Vert_{\frac{n+4}{n+2},\Omega \times (\tau,T_m)} \le K_1
\end{flalign}
Moreover, by Lemma 4.1 from \cite{morgan1989global}, we can obtain the corresponding time–trace estimate
\begin{flalign} \label{eq:time-trace}
    \|v(\tau)\|_{\frac{n+4}{n+2},(\Omega)} \le K_1,
\end{flalign}
\noindent
We claim that if $r \in \mathbb{N}$ then there exists $K_r>0$ such that 
\begin{equation}\label{eq:2.7}
\|v(\cdot,t)\|_{\left(\frac{n+4}{n+2}\right)^r,\Omega\times(0,T_m)}\le K_r
\end{equation}
To this end, note that if $1<q<\infty$ such that $q>\frac{n+2}{2}$ then $W_q^{(2,1)}(\Omega\times(\tau,T))$ embeds continuously into $C(\overline\Omega\times[\tau,T])$ and if $q<\frac{n+2}{2}$ then $W_q^{(2,1)}(\Omega\times(\tau,T))$ embeds continuously into $L^{\alpha}(\Omega\times(\tau,T))$ for all $1\le \alpha \le \frac{(n+2)q}{n+2-2q}$ \cite{ladyzhenskaia1968linear}.  Now, suppose $r\in\mathbb{N}$ such that (\ref{eq:2.7}) is true.\\
\\
\noindent
Now, we want to obtain an estimate as in (\ref{eq:2.6}) with $(\frac{n+4}{n+2})^r$ replaced by $(\frac{n+4}{n+2})^{r+1}$. Let $Q=\left(\frac{n+4}{n+2}\right)^r$, $R=\left(\frac{n+4}{n+2}\right)^{r+1}$ and $q=\frac{R}{R-1}$. If $q\ge\frac{n+2}{2}$ and $\varphi$ solves the left system in (\ref{eq:2.4}), then from above $\|\varphi\|_{p,\Omega\times(\tau,T)}$ is bounded independent of $\theta$ for every $1\le p<\infty$.  On the other hand, if $q<\frac{n+2}{2}$ then $\frac{n+4}{n+2} \le \frac{n+2+2R}{n+2}$ implies that $1 \le \frac{Q}{Q-1}\le\frac{(n+2)q}{n+2-2q}$. As a result, in either case, $\|\varphi\|_{\frac{Q}{Q-1},\Omega\times(\tau,T)}$ is bounded independent of $\theta$. So, calculating again $A,B$ and replacing them in (\ref{eq:2.5}), we can obtain bound as in (\ref{eq:2.7}) with $r$ replaced by $r+1$. Therefore, (\ref{eq:2.7}) is true. Similar argument is also applicable with (\ref{eq:time-trace}).\\
Now, let us recall that
\begin{flalign*}
    g_2(x,t) = u(x,t)v(x,t)^2-(f+\kappa)v(x,t) \le u(x,t)v(x,t)^2
\end{flalign*}
Now let us consider that $G(x,t) = u(x,t)v(x,t)^2$, then $ g_2(x,t) \le G(x,t)$. Then (\ref{eq:2.7}) implies $\Vert G(\cdot,t) \Vert_{p,\Omega \times (\tau,T_m)}$ is bounded for each $1<p<\infty$. Now, let $\Psi$ be the unique nonnegative solution to 
\begin{flalign*}
    \Psi_t(x,t) & = d_2 \Delta \Psi(x,t)+G(x,t), \ \ \ \ \ \ \Omega \times (0,T_m),\\
    \frac{\partial \Psi}{\partial \eta} & = 0, \hspace{1.4 in} \partial\Omega \times (0,T_m),\\
    \Psi & = v_0(x) \ge 0, \hspace{0.9 in} \Omega \times \{0\},
\end{flalign*}
then comparison principle for parabolic initial boundary value problems implies $v\le \Psi$, and from the $L^p(\Omega\times(0,T_m))$ bound on $G(x,t)$ for $p>\frac{n+2}{2}$,  $\|\Psi\|_{\infty,\Omega\times(0,T_m)}$ is bounded. 
Therefore,  $\Vert v(\cdot,t) \Vert_{\infty,\Omega \times (\tau,T)}$ is bounded. So, there exists a constant $\tilde{C}>0$ such that
\begin{flalign*}
\|v(\cdot,t)\|_{L^\infty(\Omega)} \le \tilde{C} \qquad \text{for all } t>0,    
\end{flalign*}
where $\tilde{C}$ depends only on $\|u_0\|_{L^\infty(\Omega)}$, $\|v_0\|_{L^1(\Omega)}$, the reaction parameters $f,k$, diffusion coefficients $d_1, d_2$, and the domain $\Omega$, but is independent of time.\\
\noindent
This along with (\ref{eq:2.2}) proves the theorem.\hfill $\square$

\section{Pattern Formation}\label{Pattern Formation}
\noindent
Let us consider either the regular Gray-Scott model (\ref{eq:local}) or mixed model (\ref{eq:1.1}) without the spatial diffusion. That is, the ODE model,
\vspace{-0.1 in}
\begin{flalign}\label{eq:3.3}
\begin{cases}
    u_t & = g_1(u,v) = -uv^2+f(1-u) \\
    v_t & = g_2(u,v) = u^2-(f+\kappa)v 
\end{cases}
\end{flalign}
By solving the steady state system i.e. $u_t=v_t=0$ we get the solutions $(u,v)=(1,0)$ and 
\vspace{-0.05 in}
\begin{flalign*}
    u_{1,2} \;=\; \frac{1}{2}\Bigl( 1 \,\pm\, \sqrt{\,1 - \frac{4\,(F+k)^{2}}{F}\,} \Bigr),
\qquad
v_{1,2} \;=\; \frac{F}{2\,(F+k)}\Bigl( 1 \,\mp\, \sqrt{\,1 - \frac{4\,(F+k)^{2}}{F}\,} \Bigr)
\end{flalign*} 
We will not go into the technical details to discuss the bifurcation analysis and condition to generate patterns by local Gray-Scott model. We refer the excellent work of \cite{mcgough2004pattern}. Pattern formation of a nonlocal reaction-diffusion system can be found in \cite{gourley2001spatio}. For fractional reaction--diffusion systems, the effects of fractional diffusion on Turing instability and pattern formation have been investigated in \cite{langlands2007turing,gafiychuk2006pattern}.\\
\noindent
For our numerical investigation, we consider, $\Omega=[0,1]^2$ with $f=0.04$, $\kappa=0.0636$, $d_1=1.0$, $d_2=0.5$, $u_0(x)=1$ and $v_0(x)=0$. We discretize the whole domain in $256 \times 256$ meshes.
We have used Finite Difference Method to simulate both local Gray-Scott system (\ref{eq:local}) and the mixed one (\ref{eq:1.1}) in for $t=100000$ sec each. Forward Euler Method has been used to approximate the integral of the kernel of fractional Laplacian operator in (\ref{eq:1.2}). The discretization of the kernal is done by a finite weighted sum over all grid points in $\Omega$, yielding a discrete nonlocal operator of convolution type. Near the boundary, the nonlocal Neumann condition is incorporated implicitly through the definition of the operator by restricting the interaction to points inside the domain, so that no additional boundary treatment is required. Central Difference Method is used for the Laplacian of the local part with zero flux ($\nabla \cdot u=\nabla \cdot v=0$) boundary conditions.\\
\noindent
Time integration is performed using a forward Euler scheme. The time step $\Delta t$ is chosen sufficiently small to satisfy the explicit stability constraints imposed by both the classical and fractional diffusion terms, which scale as $\Delta t \lesssim C h^2$ and $\Delta t \lesssim C h^{2s}$, respectively.
Stability of the scheme was verified numerically by refining $\Delta t$ and confirming that the computed patterns remained qualitatively unchanged.
\begin{figure}[H]
\centering
     \includegraphics[width=6cm]{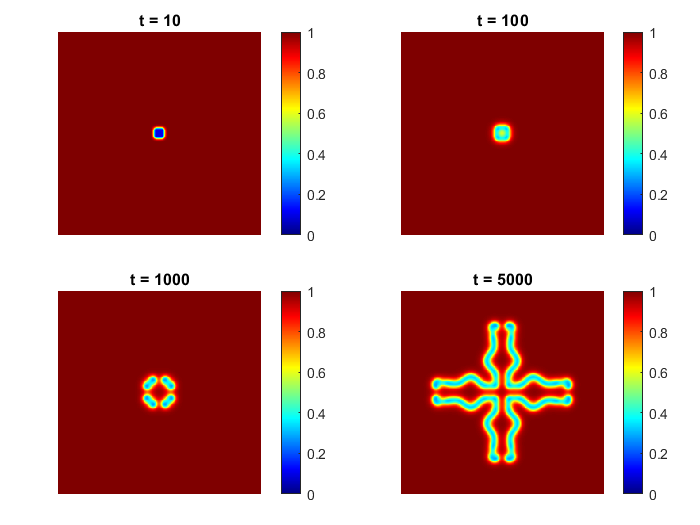}
     \includegraphics[width=6cm]{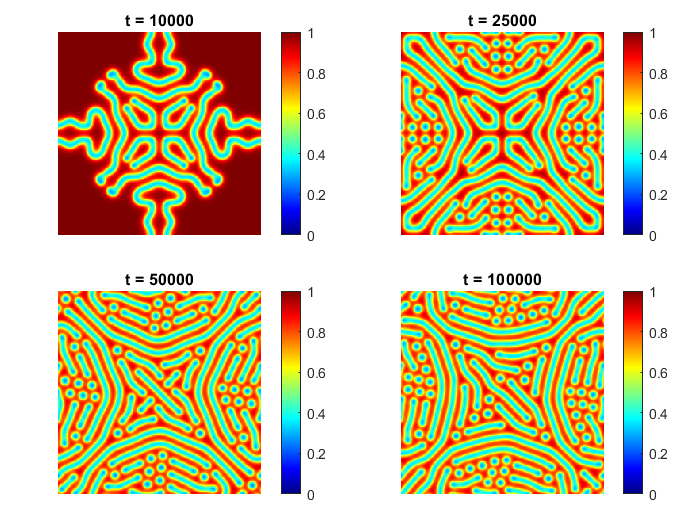}
     \caption{Pattern of local Gray-Scott model.}
    \label{fig1: Pattern of local Gray-Scott model.}
\end{figure}

\begin{figure}[H]
\centering
     \includegraphics[width=6cm]{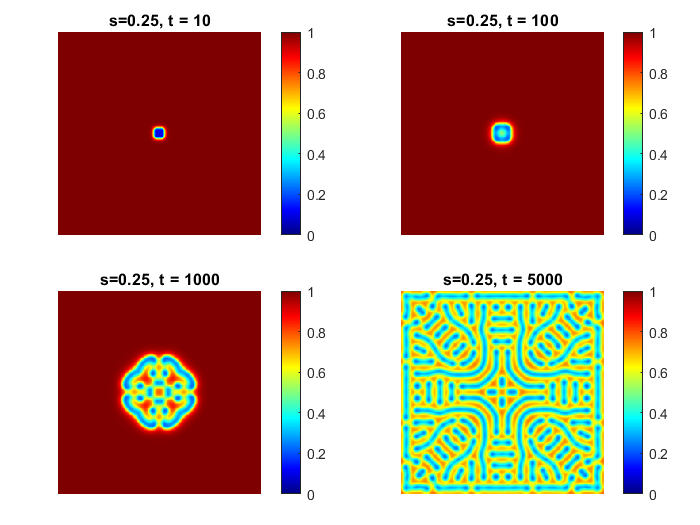}
     \includegraphics[width=6cm]{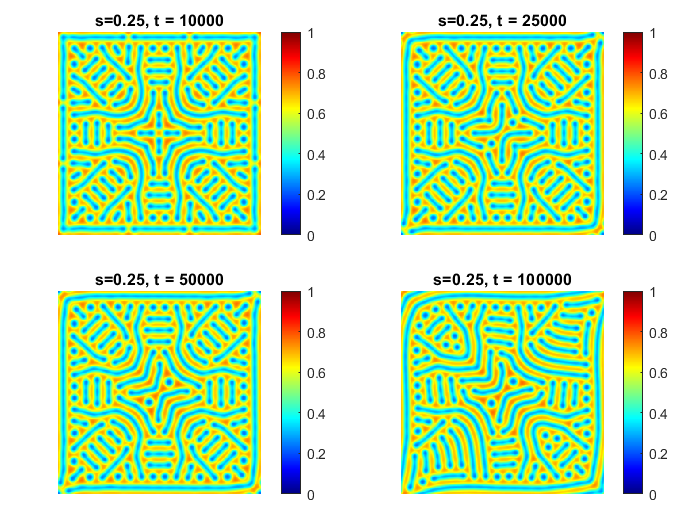}
     \caption{Pattern of mixed Gray-Scott model with $s=0.25$.}
    \label{fig2: Pattern of mixed Gray-Scott model with $s=0.25$.}
\end{figure}

\begin{figure}[H]
\centering
     \includegraphics[width=6cm]{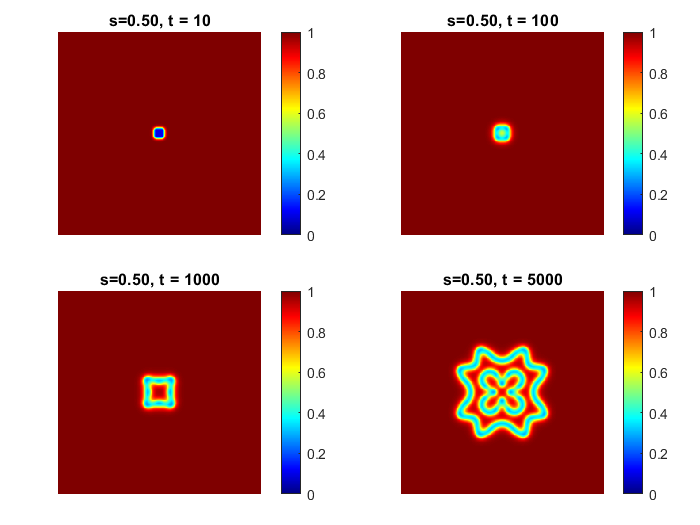}
     \includegraphics[width=6cm]{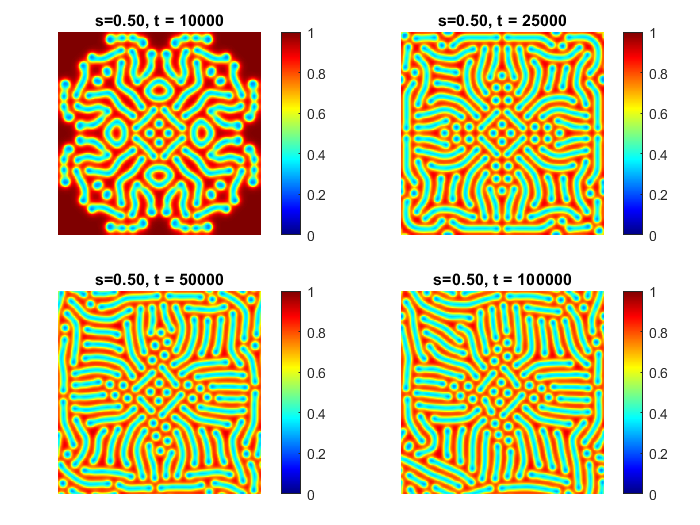}
     \caption{Pattern of mixed Gray-Scott model with $s=0.50$.}
    \label{fig3: Pattern of mixed Gray-Scott model with $s=0.50$.}
\end{figure}

\begin{figure}[H]
\centering
     \includegraphics[width=6cm]{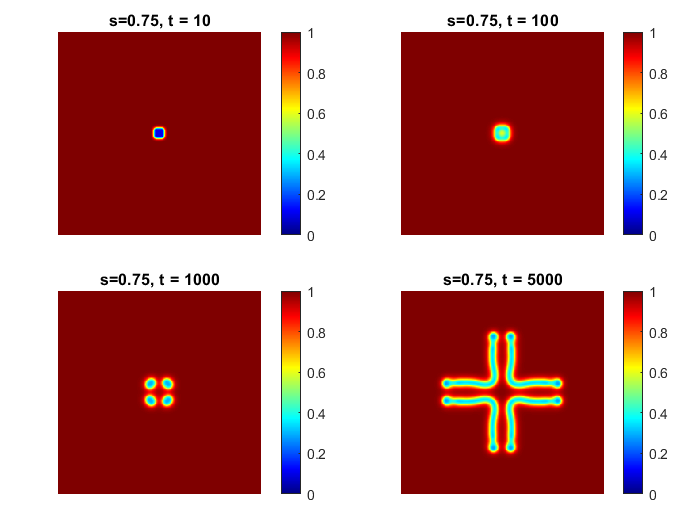}
     \includegraphics[width=6cm]{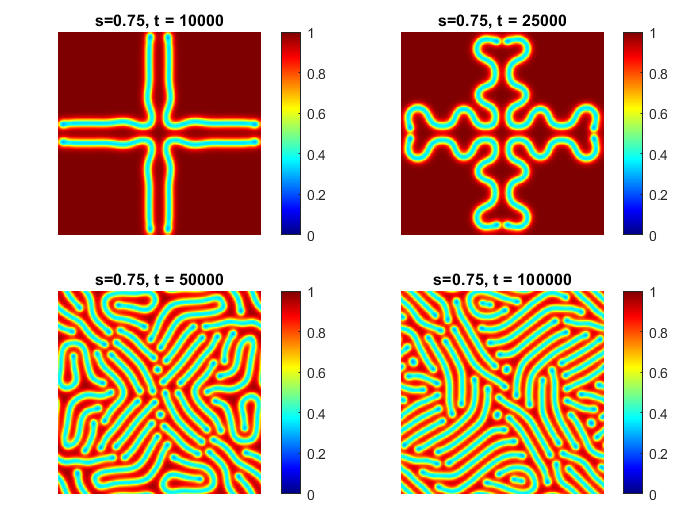}
     \caption{Pattern of mixed Gray-Scott model with $s=0.75$.}
    \label{fig4: Pattern of mixed Gray-Scott model with $s=0.75$.}
\end{figure}
\noindent
From figure (\ref{fig3: Pattern of mixed Gray-Scott model with $s=0.50$.}) and figure (\ref{fig4: Pattern of mixed Gray-Scott model with $s=0.75$.}), we observe that as $s$ increases towards $1$ the patterns starts looking similar to figure (\ref{fig1: Pattern of local Gray-Scott model.}). This numerically agrees with the analytical fact that as $s \to 1$, the fractional Laplacian coincides with regular or local Laplacian \cite{galfractional,di2012hitchhiker,bourgain2001another}. Similar results have also been observed for other nonlocal operator such as convolution operator \cite{laurencot2023nonlocal}.\\
\noindent
In our nonlocal implementation, the parameter $s$ controls the decay of the interaction kernel and hence the effective range of the diffusion.
Increasing $s$ makes the operator more local, which strengthens the damping of short spatial scales and shifts the dominant mode selection toward higher wavenumbers.
Consequently, the characteristic pattern wavelength decreases as $s$ increases; in particular, compared to figure (\ref{fig2: Pattern of mixed Gray-Scott model with $s=0.25$.}), figure (\ref{fig3: Pattern of mixed Gray-Scott model with $s=0.50$.}) and figure (\ref{fig4: Pattern of mixed Gray-Scott model with $s=0.75$.}) typically produce finer patterns with smaller dominant wavelengths, with $s=0.75$ being closer to the classical ($s\to 1$) Gray--Scott regime.

\section{Conclusion}\label{Conclusion and Future Directions}
\noindent
In this work, we study the Gray–Scott system with \emph{mixed diffusion}—a fractional Laplacian for one component and the classical Laplacian for the other—under homogeneous Neumann conditions. Using semigroup and duality techniques, we prove global existence of component-wise nonnegative solutions. Numerical experiments show that fractional diffusion qualitatively reshapes pattern morphology relative to the purely local case. These results advance the analysis of local–nonlocal models and motivate future work on stability, bifurcation, and extensions to broader nonlinearities and domains.\\
\noindent
Our analysis is carried out on a bounded smooth domain under homogeneous Neumann (zero--flux) boundary conditions, which are natural in the modeling context and allow for the standard integration-by-parts and mass-balance estimates used throughout the proofs.
The global existence and boundedness arguments extend, with only minor modifications, to other \emph{dissipative} boundary conditions for the local component (e.g.\ homogeneous Dirichlet or Robin), since the corresponding Laplacian remains sectorial and generates an analytic semigroup, and the parabolic regularity and $L^p$--bootstrap framework remain available.
However, some quantitative estimates (in particular those relying on conservation or control of spatial averages) may change, as Dirichlet conditions do not preserve mass and introduce boundary dissipation.
For the nonlocal (fractional) diffusion component, the appropriate notion of boundary condition depends on the chosen realization of the operator (restricted/regional, spectral, or censored/Neumann-type realizations), and the proofs must be adapted accordingly; in many cases, analogous semigroup and energy methods apply, but the precise form of the boundary contribution and functional setting can differ.
Finally, on unbounded domains (e.g.\ $\mathbb{R}^n$) one may still formulate the system in $L^p$-based spaces and use semigroup methods, but additional assumptions are typically required (e.g.\ on decay at infinity or integrability of initial data) to replace compactness and Poincar\'e-type inequalities that are available on bounded domains.
A detailed treatment of these extensions is beyond the scope of the present work.\\
\\
\noindent
\textbf{Conflict of Interest:} There is no conflict of Interest.\\
\\
\noindent
\textbf{Acknowledgment:} I am sincerely grateful to my advisor, Dr. Jeff Morgan, for his invaluable guidance, insightful feedback, and continued support throughout the completion of this work.
     
\bibliographystyle{acm}
\bibliography{reff}

\end{document}